\documentclass{article}%
\usepackage{amsmath}
\usepackage{amsfonts}
\usepackage{amssymb}
\usepackage{graphicx}
\usepackage{caption}
\usepackage{subfig}%
\providecommand{\U}[1]{\protect\rule{.1in}{.1in}}
\newtheorem{theorem}{Theorem}

\newtheorem{corollary}[theorem]{Corollary}

\newtheorem{lemma}[theorem]{Lemma}

\begin{document}

\begin{center}
\mbox{ }\vspace{0.5in}

{\LARGE Checking the Model and the Prior for\medskip}

{\LARGE the Constrained Multinomial}{\Large \vspace{0.25in}}

{\Large Berthold-Georg Englert\footnotemark[1], Michael Evans\footnotemark[2],
Gun Ho Jang\footnotemark[3], Hui Khoon Ng\footnotemark[4], David
Nott\footnotemark[5] and Yi-Lin Seah\footnotemark[6]\bigskip}
\end{center}

\noindent\textbf{Abstract:} The multinomial model is one of the simplest
statistical models. When constraints are placed on the possible values for the
probabilities, however, it becomes much more difficult to deal with. Model
checking and checking for prior-data conflict is considered here for such
models. A theorem is proved that establishes the consistency of the check on
the prior. Applications are presented to models that arise in quantum state
estimation as well as the Bayesian analysis of models for ordered probabilities.


\footnotetext[1]{Centre for Quantum Technologies (CQT), Dept. of Physics,
National University of Singapore and MajuLab, Singapore}\footnotetext[2]{Dept.
of Statistical Sciences, University of Toronto, Toronto, Ontario,\ Canada}%
\footnotetext[3]{Ontario Institute for Cancer Research, Toronto, Ontario,
Canada}\footnotetext[4]{Yale-NUS College, CQT and MajuLab, Singapore}%
\footnotetext[5]{Dept. of Statistics and Applied Probability, National
Unversity of Singapore, Singapore}\footnotetext[6]{Centre for Quantum
Technologies, Singapore}

\section{Introduction}

Suppose we have a sample of $n$ from a multinomial$(1,\theta_{1},\ldots
,\theta_{k+1})$ distribution where $\theta=(\theta_{1},\ldots,\theta_{k}%
)\in\Theta_{k}$ is an unknown element of
\[
\Theta_{k}=\{(\theta_{1},\ldots,\theta_{k}):\theta_{1}+\cdots+\theta
_{k}<1,\theta_{i}>0,i=1,\ldots,k\}
\]
and $\theta_{k+1}=1-\theta_{1}-\cdots-\theta_{k}$. Note that throughout the
paper the notation $\theta$\ always refers to first $k$ probabilities. If
$T_{n}=(T_{1},\ldots,T_{k})$ denotes the counts from the first $k$ categories,
with $T_{k+1}=n-T_{1}-\cdots-T_{k}$, then $T_{n}$ is a minimal sufficient
statistic (mss) and $T_{n}\sim\,$multinomial$(n,\theta_{1},\ldots,\theta
_{k+1}).$

In some applications something is known about the true value of $\theta,$
beyond the fact that it is in $\Theta_{k},$ and this is expressed in the form
of a prior probability distribution $\Pi$ on $\Theta_{k}.$ It is assumed that
the prior is absolutely continuous with respect to volume measure on
$\Theta_{k}$ with the density of $\Pi$ denoted by $\pi.$ The prior can be
thought of as a way of imposing constraints on $\theta.$ These may be soft
constraints such that $\pi(\theta)$ is relatively low in parts of $\Theta
_{k},$ reflecting beliefs that these values are improbable, or hard
constraints where $\pi(\theta)=0$ for values that are `known' to be
impossible. The prior is to be thought of as the marginal distribution of
$\theta$ and the multinomial$(n,\theta_{1},\ldots,\theta_{k+1})$ is the
conditional for $T_{n}$ given $\theta,$ together giving rise to the joint
probability model for $(\theta,T_{n}).$

Given that the model and prior are appropriate in a particular application,
the inferential analysis proceeds by first obtaining the conditional
distribution for $\theta$ given $T_{n},$ the posterior of $\theta.$ The
posterior represents beliefs about the true value of $\theta$ after observing
the data. While a variety of approaches can be considered for inference about
$\theta,$ our concern here is with whether or not the specified model and
prior are appropriate.

The primary way to determine whether or not an ingredient to a statistical
analysis is appropriate is to compare it somehow with the data. For example,
if the observed data is surprising for every distribution in the model, then
it is reasonable to question the appropriateness of the model. By surprising
it is meant that the data falls in a region such that each distribution in the
model gives a relatively low probability for that data's occurrence. This
assessment is usually approached via the computation of a p-value in a
so-called goodness-of-fit test. For example, for the multinomial this could be
assessed by a generalization of the well-known runs test as this is assessing
whether or not the unreduced data is i.i.d. In this paper it will always be
assumed that i.i.d. sampling holds. This can be induced (approximately) when
random sampling from large populations but otherwise needs to be checked.

It is to be noted that model checking procedures are typically based on
aspects of the data beyond the values of the mss although, as subsequently
discussed, some qualifications are necessary. So, for example, the
goodness-of-fit test could be based on the conditional distribution of the
original data given the value of the minimal sufficient statistic. This
conditional distribution is completely independent of $\theta$ and so data
that was surprising for this distribution is indicating a problem with the
model. In certain cases there are ancillaries that are functions of the mss,
however, and then model checking can proceed by comparing the values of these
ancillaries with their known distributions. For example, the introduction of
hard-constraints with the multinomial can indeed produce such ancillaries.

In general, there are many ways in which model checking can proceed and it
doesn't seem possible to argue definitively for one approach over another.
There are, however, some basic principles that seem necessary. For example,
being careful about what aspects of the data are used in model checking seems
paramount. Another basic principle would seem to be the separation of the
checking of the prior from checking the model. If the ingredients fail the
tests of appropriateness, then we would like to know specifically what
component caused the failure. If we try to jointly check the model and prior,
say by some aspect of the posterior, then failure cannot be assigned to the
specific component. Also, given that the prior is implicitly dependent on the
model, it isn't meaningful to check the prior if the model fails its checks.
So as argued in Evans and Moshonov (2006), it makes sense to separate the
checks on the model and prior and perform the check on the model first. Also,
it is shown in Evans and Moshonov (2006, 2007) that the check on the prior can
sometimes be decomposed so that individual components of the prior, as when
the prior is specified hierarchically, can be checked so that an isolated
aspect of the prior can be identified as causing the problem when one exists.

While the separation of the check on the prior and the model is a principle
worth noting, this does not rule out the possibility of a Bayesian approach to
model checking. For example, suppose hard constraints lead to the model being
a subclass $M^{\prime}$ of the full multinomial family $M$. It is then
reasonable to place a uniform prior on $M,$ equivalently a uniform prior on
$\theta,$ so that each multinomial has the same weight and then assess whether
or not $M^{\prime}$ is `plausible'. Just how this assessment is to be carried
out is a matter for some discussion but our preference is a comparison of the
prior belief in $M^{\prime}$ with its posterior belief. So if belief in
$M^{\prime}$ has increased after seeing the data, then there is evidence for
$M^{\prime}$ being true and if it has decreased, then there is evidence
against $M^{\prime}$ being true.\ Such an approach is discussed in Al-Labadi
and Evans (2016) and Al-Labadi, Baskurt and Evans (2017) where a distance
measure is introduced and the concentration of the posterior about $M^{\prime
}$ is compared with the concentration of the prior about $M^{\prime}$ to
assess whether or not there is evidence for or against $M^{\prime}.$ A notable
aspect of this, is that while the check is Bayesian and based on a measure of
statistical evidence, it does not involve the prior $\Pi$ which is only
checked if there is evidence in favor of $M^{\prime}.$ Note too that this
approach involves conditioning on all the data and so avoids the issues that
arise when there are multiple maximal ancillaries.

In Section 2 model checking is discussed for the constrained multinomial,
namely, when $\theta\in\Theta$ and $\Theta$ is a proper subset of $\Theta
_{k}.$ Constrained multinomial models arise naturally in quantum state
estimation and such an example is developed in the paper as well as
goodness-of-fit for the multinomial model with ordered probabilities. In
Section 3 methods are discussed for checking the prior and a consistency
result is established for the specific check used that substantially
generalizes a result established in Evans and Jang (2011b). Also, an
elicitation algorithm is developed for the multinomial with ordered
probabilities and this leads to a methodology for checking elicited
information which is somewhat different than checking a specific prior.

\section{Checking the Model}

Suppose that we are satisfied that the data is i.i.d. multinomial$(1,\theta
_{1},\ldots,\theta_{k+1})$ for some $\theta=(\theta_{1},\ldots,\theta_{k}%
)\in\Theta_{k}.$ Suppose further, however, it is believed that $\theta
^{\text{true}}\in\Theta,$ where $\Theta$ is a proper subset of $\Theta_{k},$
and it is desirable for the analysis to reflect this. It is then necessary to
check that indeed this constraint is appropriate for the data obtained. For
example, it might be that the data was not collected correctly and this led to
a failure of the model such that, while the multinomial assumption is correct,
$\theta^{\text{true}}\notin\Theta$. Note that for the model in question any
elicited prior $\Pi$ on $\theta$ must satisfy $\Pi(\Theta)=1.$ By
\textit{elicited prior }is meant a prior that is selected based upon a process
that reflects what is known about the true distribution being sampled from.

Consider the following examples.\smallskip

\noindent\textbf{Example 1. }\textit{Quantum state estimation. }

Counts of events are recorded associated with the state of a qubit. A quantum
measurement involves $k+1$ detectors and quantum theory leads to a
distribution that corresponds to sampling from a multinomial$(1,\theta
_{1},\ldots,\theta_{k+1})$ for some $k$ where $\theta_{i}$ is the probability
that a \textquotedblleft click\textquotedblright\ is recorded from the $i$-th
detector$.$ Quantum theory dictates, depending on how the measurements are
taken, that the probabilities satisfy certain constraints.

For example, when $k+1=3,$ then the \textit{symmetric trine} model imposes the
constraint $\theta_{1}^{2}+\theta_{2}^{2}+\theta_{3}^{2}\leq1/2,$ with more
involved constraints required for the asymmetric case as discussed in Example
3. When $k+1=4,$ then the \textit{cross-hairs} model imposes the constraints
$\theta_{1}+\theta_{2}=1/2,\theta_{3}+\theta_{4}=1/2$ and $\theta_{1}%
^{2}+\theta_{2}^{2}+\theta_{3}^{2}+\theta_{4}^{2}\leq3/8,$ while the
\textit{tetrahedron} model corresponds to the constraint $\theta_{1}%
^{2}+\theta_{2}^{2}+\theta_{3}^{2}+\theta_{4}^{2}\leq1/3$ only. When $k+1=6,$
then the \textit{Pauli} model imposes the constraints $\theta_{1}+\theta
_{2}=1/3,\theta_{3}+\theta_{4}=1/3,\theta_{5}+\theta_{6}=1/3$ and $\theta
_{1}^{2}+\theta_{2}^{2}+\theta_{3}^{2}+\theta_{4}^{2}+\theta_{5}^{2}%
+\theta_{6}^{2}\leq2/9$. More on how these models arise can be found in Shang,
Ng, Sehrawat, Li and Englert (2013) but, sufficed to say, these applications
produce a rich variety of constrained multinomial models.\smallskip

\noindent\textbf{Example 2. }\textit{Contingency tables with ordered
probabilities.}

In some circumstances it is reasonable to suppose that the probabilities
satisfy an ordering such as
\begin{equation}
\theta_{1}\geq\theta_{2}\geq\cdots\geq\theta_{k+1}. \label{ordered}%
\end{equation}
Such a model arises in contexts where systems exhibit aging as in, for
example, Briegel, Englert, Sterpi and Walther (1994). Issues associated with
checking (\ref{ordered}), and with eliciting a prior on $\theta$ when this
restriction is deemed correct, lead to the utility of the following
result.\smallskip

\noindent\textbf{Lemma 1}. Any $\theta\in\Theta_{k}$ satisfying (\ref{ordered}%
) is given by, for some $\omega\in\Theta_{k},$%
\begin{equation}
(\theta_{1},\ldots,\theta_{k+1})^{t}=A_{k}(\omega_{1},\ldots,\omega_{k+1})^{t}
\label{basiceq}%
\end{equation}
and any $\omega\in\Theta_{k}$ produces a $\theta\in\Theta_{k}$ satisfying
(\ref{ordered}) via (\ref{basiceq}), with%
\[
A_{k}=\left(
\begin{array}
[c]{ccccc}%
1 & 1/2 & 1/3 & \ldots & 1/(k+1)\\
0 & 1/2 & 1/3 & \ldots & 1/(k+1)\\
0 & 0 & 1/3 & \ldots & 1/(k+1)\\
\vdots & \vdots & \vdots & \vdots & \vdots\\
0 & 0 & 0 & \ldots & 1/(k+1)
\end{array}
\right)  .
\]

\noindent\textbf{Proof}: First assume $(\theta_{1},\ldots,\theta_{k+1}%
)^{t}=A_{k}(\omega_{1},\ldots,\omega_{k+1})^{t}$ for $\omega\in\Theta_{k}.$
Then $\theta_{i}=\sum_{j=i}^{k+1}\omega_{j}/j$ and it is clear that
$0\leq\theta_{i}\leq1$ and $\theta_{1}+\cdots+\theta_{k+1}=\omega_{1}%
+\cdots+\omega_{k+1}=1$ so $\theta\in\Theta_{k}.$ It is also immediate that
$\theta$ satisfies (\ref{ordered}).

Now suppose $\theta\in\Theta_{k}$ satisfies (\ref{ordered}) and put
$(\omega_{1},\ldots,\omega_{k+1})^{t}=A_{k}^{-1}(\theta_{1},\ldots
,\theta_{k+1})^{t}.$ An easy calculation shows that
\[
A_{k}^{-1}=\left(
\begin{array}
[c]{ccccc}%
1 & -1 & 0 & \ldots & 0\\
0 & 2 & -2 & \ldots & 0\\
0 & 0 & 3 & \ldots & 0\\
\vdots & \vdots & \vdots & \vdots & \vdots\\
0 & 0 & 0 & \ldots & k+1
\end{array}
\right)  .
\]
Therefore, $\omega_{i}=i(\theta_{i}-\theta_{i+1})$ for $i=1,\ldots,k$ and
$\omega_{k+1}=(k+1)\theta_{k+1}.$ Since $\theta_{i}\geq\theta_{i+1},$ then
$\omega_{i}\geq0$ and if $\omega_{i}>1$ then $\theta_{i}>\theta_{i+1}+1/i>1/i$
which by (\ref{ordered}) implies $\theta_{1}+\cdots+\theta_{i}>1$ which is
false.\ So $\omega_{i}\in\lbrack0,1]$ for $i=1,\ldots,k$ and similarly
$\omega_{k+1}\in\lbrack0,1].$ Finally, $\omega_{1}+\cdots+\omega_{k+1}%
=\sum_{j=1}^{k}i(\theta_{i}-\theta_{i+1})+(k+1)\theta_{k+1}=\theta_{1}%
+\cdots+\theta_{k+1}=1$ and this completes the proof.\smallskip\ 

A particular model satisfying (\ref{ordered}) is given by the Zipf-Mandelbrot
distribution where, for parameters $\alpha>-1,\beta\geq0,$%
\begin{equation}
\theta_{i}=(\alpha+i)^{-\beta}/C_{k}\left(  \alpha,\beta\right)
\label{ZMprobs}%
\end{equation}
for $i=1,\ldots,k+1$ where $C_{k}\left(  \alpha,\beta\right)  =\sum
_{i=1}^{k+1}(\alpha+i)^{-\beta}.$ When $\beta=0$ this is the uniform
distribution and for fixed $\beta$ this converges to the uniform as
$\alpha\rightarrow\infty.$ For fixed $\alpha$ the distribution becomes
degenerate on the first cell as $\beta\rightarrow\infty.$ For large
$k,\alpha=0$ and $\beta>1$ the zeta$(\beta)$ distribution, with $\theta
_{i}\propto i^{-\beta}$ for $i=1,2,\ldots$ serves as an approximation. As
discussed in Izs\'{a}k (2006), there are a variety of applications of this
distribution as in word frequency distributions in texts. The distribution
given by (\ref{ZMprobs}) is denoted here as ZM$_{k}\left(  \alpha
,\beta\right)  .$

In dose-response models, as discussed in Chuang-Stein and Agresti (1997),
there are $I$ treatments, possibly corresponding to the frequency of a
particular treatment, and $J$ response classifications reflecting the severity
of a condition from nonexistent to most severe. When there is a monotone
increasing effect based on say frequency of treatment, it makes sense to
assume $\theta_{1j}\leq\theta_{2j}\leq\cdots\leq\theta_{k+1j}$ for
$j=1,\ldots,J.$\smallskip

\noindent So there are many applications of the constrained multinomial.

Following the discussion in Section 1 the constrained multinomial model is
checked first and, if the model passes, the prior $\Pi$ placed on the
restricted parameter space $\Theta$ is then checked. The check on the
constraints should not involve $\Pi$ in any way since it is not involved in
the production of the data. One of the advantages of this is that we can
contemplate using another prior $\Pi_{0}$ on $\Theta_{k}$ where supp$(\Pi
_{0})=\Theta_{k},$ for the model checking step. This leads to formal inference
methods for model checking. For example, taking $\Pi_{0}$ to be the uniform
prior on $\Theta_{k}$ is quite natural as this treats all multinomials
equivalently and so this prior is used for the model checking step hereafter.

Suppose first that $\Pi_{0}(\Theta)>0.$ Based on the observed $T_{n}%
=T_{n}(x),$ the posterior probability of $\Theta$ is
\[
\Pi_{0}(\Theta\,|\,T_{n}(x))=\int_{\Theta}f_{\theta}(T_{n}(x))\,\pi_{0}%
(\theta)\,d\theta/m_{T_{n}}(T_{n}(x))
\]
where $f_{\mathbf{\theta}}$ is the multinomial$(n,\theta_{1},\ldots
,\theta_{k+1})$ density, $\pi_{0}$ is the Dirichlet$(1,\ldots,1)$ density and
$m_{T_{n}}(T_{n}(x))=\int_{\Theta_{k}}f_{\mathbf{\theta}}(T_{n}(x))\,\pi
_{0}(\theta)\,d\theta.$ Note that the posterior distribution of $\theta$ is
Dirichlet$(T_{1}(x)+1,\ldots,T_{k+1}(x)+1)$ which can be used in evaluating
$\Pi_{0}(\Theta\,|\,T_{n}(x)).$ The relative belief ratio in favor of the
hypothesis $H_{0}=\Theta,$ namely, the hypothesis that the constraints hold,
is then
\[
\text{RB}(\Theta\,|\,T_{n}(x))=\Pi_{0}(\Theta\,|\,T_{n}(x))/\Pi_{0}(\Theta)
\]
and there is evidence in favor of $H_{0}$ when RB$(\Theta\,|\,T_{n}(x))>1$ and
evidence against when RB$(\Theta\,|\,T_{n}(x))<1.$ It seems reasonable in such
circumstances to not be concerned about the model if RB$(\Theta\,|\,T_{n}%
(x))>1$ and consider that a problem has occurred whenever RB$(\Theta
\,|\,T_{n}(x))<1.$ It is possible, however, to also assess the strength of
this evidence by reporting $\Pi_{0}(\Theta\,|\,T_{n}(x))$ as this represents
our belief (as opposed to the evidence) that $H_{0}$ is true. So, if
RB$(\Theta\,|\,T_{n}(x))>1$ and $\Pi_{0}(\Theta\,|\,T_{n}(x))$ is low, then
there is only weak evidence in favor of $H_{0}$ while, if RB$(\Theta
\,|\,T_{n}(x))<1$ and $\Pi_{0}(\Theta\,|\,T_{n}(x))$ is low, then there is
strong evidence against $H_{0}.$ Similarly, if RB$(\Theta\,|\,T_{n}(x))>1$ and
$\Pi_{0}(\Theta\,|\,T_{n}(x))$ is high, then there is strong evidence in favor
of $H_{0}$ while, if RB$(\Theta\,|\,T_{n}(x))<1$ and $\Pi_{0}(\Theta
\,|\,T_{n}(x))$ is high, then there is only weak evidence against $H_{0}.$
Note that it is important to separate the measurement of evidence from the
measurement of belief as is done here.

In some circumstances it may arise that $\Pi_{0}(\Theta)=0$ simply because
$\Theta$ is a lower dimensional subset of $\Theta_{k}.$ As such, we cannot
proceed as just described. In fact, if $\Pi_{0}(\Theta)$ is very small, then
again the preceding approach to checking the model seems questionable as it
cannot be expected that $\Pi_{0}(\Theta\,|\,T_{n}(x))$ will be large, and so
obtain strong evidence in favor of the model when that is appropriate, unless
the amount of data is very large. To deal with these problems consider the
approach discussed in Al-Labadi and Evans (2016) and Al-Labadi, Baskurt and
Evans (2017). For this $d_{H_{0}}:H_{0}\rightarrow\lbrack0,\infty)$ is
specified where $d_{H_{0}}(\theta)$ is a measure of the distance of
$\mathbf{\theta}$ from $H_{0}$ with $d_{H_{0}}(\theta)=0$ iff $H_{0}$ is true.
For example, in some contexts squared Euclidean distance might make sense so
that $d_{H_{0}}(\theta)=\inf_{\mathbf{\theta}^{\prime}\in\Theta}%
||\theta-\theta^{\prime}||^{2}$ for $\theta\in\Theta_{k}.$ Perhaps a more
natural measure is obtained using Kullback-Leibler divergence, namely,
$d_{H_{0}}(\theta)=\inf_{\mathbf{\theta}^{\prime}\in\Theta}$KL$(\theta
\,||\,\theta^{\prime}).$

The hypothesis $H_{0}$ is then assessed via the relative belief ratio
\[
\text{RB}_{d_{H_{0}}}(0\,|\,T_{n}(x))=\lim_{\delta\downarrow0}\text{RB}%
_{d_{H_{0}}}([0,\delta)\,|\,T_{n}(x)).
\]
Typically the limit cannot be computed exactly so $\delta>0$ is selected such
that RB$_{d_{H_{0}}}([0,\delta)\,|\,T_{n}(x))\approx\,$RB$_{d_{H_{0}}%
}(0\,|\,T_{n}(x)).$ In practice, there is a $\delta>0$ such that, if
$d_{H_{0}}(\theta_{true})\in\lbrack0,\delta),$ then $H_{0}$ can be regarded as
true. The value of $\delta$ can be determined via bounding the absolute error
in the probabilities (squared Euclidean distance) or bounding the relative
error in the probabilities (KL divergence), see Al-Labadi, Baskurt and Evans
(2017). So, if RB$_{d_{H_{0}}}([0,\delta)\,|\,T_{n}(x))>1$ there is evidence
in favor of $H_{0}$ and evidence against when RB$_{d_{H_{0}}}([0,\delta
)\,|\,T_{n}(x))<1.$ The strength of the evidence can be measured by
discretizing the range of the prior distribution of $d_{H_{0}}$ into
$[0,\delta),[\delta,2\delta),\ldots$ and then computing the posterior
probability $\Pi_{0}(\{i:\,$RB$_{d_{H_{0}}}([0,\delta)\,|\,T_{n}(x))\leq
\,$RB$_{d_{H_{0}}}([0,\delta)\,|\,T_{n}(x))\}\,|\,T_{n}(x)).$

The consistency of this approach to model checking follows from results in
Evans (2015). As $n\rightarrow\infty$ the relative belief ratio converges to
its maximum possible value (greater than 1) and the strength goes to 1 when
$H_{0}$ is true and the relative belief ratio and the strength go to 0 when
$H_{0}$ is false.

Consider now applications to some examples.\smallskip

\noindent\textbf{Example 3. }\textit{Goodness-of-fit for the trine model.}

Table \ref{datatab} contains data from two separate experiments discussed in
Len, Dai, Englert and Krivitsky (2017) where models corresponding to two
instances of the trine model are relevant and $n_{i}$ is the number of clicks
on the $i$-th detector.
\begin{table}[tbp] \centering
\begin{tabular}
[c]{|l|c|c|c|c|}\hline
& $n$ & $n_{1}$ & $n_{2}$ & $n_{2}$\\\hline
Symmetric & \multicolumn{1}{|r|}{7076} & \multicolumn{1}{|r|}{3416} &
\multicolumn{1}{|r|}{1912} & \multicolumn{1}{|r|}{1748}\\
Asymmetric & \multicolumn{1}{|r|}{6756} & \multicolumn{1}{|r|}{6192} &
\multicolumn{1}{|r|}{316} & \multicolumn{1}{|r|}{248}\\\hline
\end{tabular}
\caption{Results from two experiments based on the trine model in Example 3.}\label{datatab}%
\end{table}
For these models $\Theta=\{\theta:\left(  \theta-c\right)  ^{t}C\left(
\theta-c\right)  \leq1\}$ where%
\[
c=\frac{1}{2}\left(
\begin{array}
[c]{c}%
2a\\
1-a
\end{array}
\right)  ,C=(1-2a)^{-1}\left(
\begin{array}
[c]{cc}%
(1-1/a)^{2} & 2\\
2 & 4
\end{array}
\right)  ,
\]
$a=0.5\sin^{2}(\cos^{-1}(\cot(2\varphi_{0})))$ and $\varphi_{0}$ an angle
associated with the experiment.

For the symmetric trine $\varphi_{0}=\pi/6$ so $a=1/3$ and $\Pi_{0}%
(\Theta)=a\sqrt{1-2a}\pi=0.6046.$ Under $\Pi_{0}$ the posterior distribution
of $\theta$ is Dirichlet$(3417,1913,1749)$ and sampling from this distribution
shows that the entire posterior is concentrated within $\Theta$ so
RB$(\Theta\,|\,T_{n}(x))=1/0.6046=1.6540.$ So there is evidence in favor of
the symmetric trine model and this is very strong evidence since the posterior
content of $\Theta$ is effectively 1.

For the asymmetric trine case $\varphi_{0}=2\pi/9$ so $a=0.48445$ and $\Pi
_{0}(\Theta)=0.2684.$ Under $\Pi_{0}$ the posterior distribution of $\theta$
is Dirichlet$(6193,317,249)$ and sampling from this distribution shows that
the entire posterior is concentrated within $\Theta$ so RB$(\Theta
\,|\,T_{n}(x))=1/0.2684=3.7258.$ So there is evidence in favor of the
symmetric trine model and again this is very strong evidence.

So with both models one can feel quite confident that the true values of the
probabilities lie within the respective sets $\Theta.$ The evidence is pretty
definitive here because of the large amount of data collected.\smallskip

\noindent\textbf{Example 4. }\textit{Goodness-of-fit for ordered
probabilities.}

A numerical example used in Izs\'{a}k (2006), based on data concerned with fly
diversity found in Papp (1992), is considered where the counts are given by
$f=(145,96,35,29,20,11,4,4,4,3,3,2,2,1,1,1,1,1).$ The question is whether or
not these data can reasonably be thought of as coming from the model given by
(\ref{ordered}) and even from the submodel given by the collection of
Zipf-Mandelbrot distributions. So here $k=17$ and $n=363$, the prior $\Pi_{0}$
is $\theta\sim\,$Dirichlet$(1,\ldots,1)$ and the posterior is $\theta\sim
\,$Dirichlet$(146,97,36,\ldots,2).$

Consider first checking (1). For this model $\Pi_{0}(\Theta
)=1/18!=1.5619\times10^{-16}$ which is exceedingly small and this suggests
that estimating $\Pi_{0}(\Theta\,|\,f)$ with any accuracy will be very
difficult. It is to be noted, however, that given the small prior probability
of this set, if any of the values generated from the posterior for some
feasible sample size fall in $\Theta$, then this will give clear evidence in
favor of the model. For example, in a sample of $10^{7}$ the posterior content
was estimated as $10^{-7}$ and this produces a relative belief ratio of
$6.4\times10^{8}$ but the strength of this evidence in favor is exceedingly
weak. A better approach in such a problem is to group the cells into
sequential subgroups such that the hypothesized monotonicity in the model is
maintained. For example, here there are 18 cells and so 9 groups of size 2 or
6 groups of size 3 are possible. Also, 4 groups of size 4 are possible with a
fifth group of size 2. Clearly it is always possible to group the cells in
this way so that monotonicity in the probabilities is maintained. To select
which grouping to use for the test it makes sense to start with the finest
grouping such that the posterior content can be accurately estimated but
coarser groupings can also be examined. Choosing 9 groups of size 2 worked
here as the posterior content of the relevant set was estimated as
$0.0396,0.0402$ and $0.0406$ based on samples of sizes $10^{4},10^{5}$ and
$10^{6},$ respectively. Since the prior content of the relevant set is
$1/9!=2.755\,7\times10^{-6},$ the relative belief ratio is $14726$ although,
based on the posterior content, this appears to be only weak evidence in
favor. A coarser grouping leads to increased confidence in the model. For
example, with groups of size 3 the relative belief ratio is $285.6312$ and the
posterior content is $0.40.$

Consider now checking the ZM$_{k}$ model. Since the set of ZM$_{k}$
distributions has prior probability 0 with respect to the uniform
distribution, the distance measure approach is necessary and the KL distance
measure is used. A technical difficulty involved here is the need to find, for
given $\theta,$ the value of%
\begin{equation}
d_{H_{0}}(\theta)=\inf_{\alpha>-1,\beta\geq0}\sum_{i=1}^{k+1}\theta_{i}%
\ln(\theta_{i}C_{k}\left(  \alpha,\beta\right)  (\alpha+i)^{\beta})
\label{kl2}%
\end{equation}
for each generated $\theta,$ to obtain samples from the prior and posterior
distribution of $d_{H_{0}}.$ For this a large table of ZM$_{k}$ distributions
was created, (\ref{kl2}) computed for each element and the distribution found
that minimizes this quantity.

As discussed in Example 2, there is redundancy in the parameterization of this
family. For example, uniformity is well-approximated by many values of
$\left(  \alpha,\beta\right)  .$ A value of $\delta>0$ is selected, however,
so that if the KL distance between two distributions is less than $\delta,$
then this difference is irrelevant for the application. Note that $\log
(\theta_{i}/p_{i})=\log(1+(\theta_{i}-p_{i})/p_{i})\approx(\theta_{i}%
-p_{i})/p_{i}$ when $(\theta_{i}-p_{i})/p_{i}$ is small and so, if
KL$(\theta,p)=\sum_{i=1}^{k+1}\theta_{i}\log(\theta_{i}/p_{i})\leq\delta,$
then the average relative difference between the probabilities given by
$\theta\ $and $p$ is immaterial. So for any $\beta\geq0$ it is only necessary
to consider values of $\alpha$ such that the KL distance
\begin{equation}
(k+1)^{-1}\sum_{i=1}^{k+1}\ln(C_{k}\left(  \alpha,\beta\right)  (\alpha
+i)^{\beta}/(k+1)) \label{kldist}%
\end{equation}
is greater than or equal to $\delta$\ and this places an upper bound on
$\alpha.$ As such, this redundancy plays no role in assessing the goodness-of-fit.

Figure 1 is plot of the prior and posterior densities of $d_{H_{0}}$ based on
Monte Carlo samples of size $10^{5},$ using $\delta=0.02$ and some smoothing.
It is clear from this that the posterior has become much more concentrated
about the ZM$_{k}$ model than the prior. Furthermore, RB$(([0,\delta
)\,|\,f)=1.75\times10^{3}$ and the strength equals 1. So there is ample
evidence in favor of the ZM$_{k}$ model with this data set. This agrees
somewhat with the finding in Izs\'{a}k (2006) who conducted a goodness-of-fit
test via computing a p-value based on the chi-squared statistic after grouping
and found no evidence against the model. Note that with the methodology
developed here, there is no need to appeal to asymptotics.

It may seem anomalous that stronger evidence is found for the ZM$_{k}$ model
than for the bigger model of ordered probabilities. One might try to account
for this by noting that different methodologies are used in the two cases, but
even when the same methods are used it is possible to have evidence for a
subset and evidence against the superset let alone just weaker evidence. This
phenomenon can arise, when the prior probability of the subset is relatively
small when compared to the prior probability of the superset as is the case
here. So evidence in favor of a subset may be mitigated by the other
possibilities in the bigger set. A full discussion and relevant example
concerning this can be found in Evans (2015).%
\begin{figure}
[ptb]
\begin{center}
\includegraphics[
height=2.4708in,
width=2.4708in
]%
{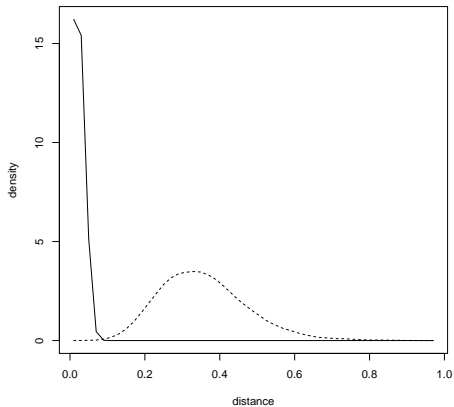}%
\caption{Plot of prior (- - -) and posterior (---) densities of $d_{H_{0}}$ in
Example 4.}%
\end{center}
\end{figure}

\section{Checking the Prior}

Suppose the model has been found to be acceptable so that now attention
focuses on the elicited prior $\Pi.$ A prior is inappropriate when the prior
places relatively little mass in a region containing the true value. Of course
the true value is not known but still there are a number of ways in which the
prior can be checked. Overall, however, all methods for checking the prior are
assessing whether or not such a contradiction exists.

The approach taken in Evans and Moshonov (2006) is used here. The prior
predictive density of $T_{n}$ is given by
\begin{equation}
m_{T_{n}}(\mathbf{t})=\frac{n!}{\prod_{j=1}^{k+1}t_{j}!}\int_{\Theta_{k}}%
\prod_{j=1}^{k+1}\theta_{j}^{t_{j}}\,\pi(\theta)\,d\theta. \label{dens}%
\end{equation}
Based on this density the check on the prior is to compute
\begin{equation}
M_{T_{n}}(m_{T_{n}}(\mathbf{t})\leq m_{T_{n}}(T_{n}(x))). \label{pvalue}%
\end{equation}
Clearly (\ref{pvalue}) is measuring where the observed $T_{n}(x)$ lies with
respect to its prior distribution as $m_{T_{n}}$ is the prior density of
$T_{n}.$ If (\ref{pvalue}) is small, then $T_{n}(x)$ lies in a region of low
prior probability for $T_{n}$ and it is apparent that the data and the prior
are in conflict. It is also clear that the check on the prior should depend on
the data only through a mss because the conditional distribution of the
remaining aspects of the data beyond the mss does not depend on $\theta$ and
so can reveal nothing about the adequacy of the prior.

In Evans and Jang (2011b) a general consistency result is established for
(\ref{pvalue}) as it is shown there that, under conditions, this converges to
$\Pi(\pi(\theta)\leq\pi(\theta^{\text{true}})).$ A small value of
(\ref{pvalue}) is thus an indication that the prior is placing its mass in the
wrong place as this suggests that the true value lies in a region of
relatively low prior probability. One of the conditions for convergence,
however, is that the prior predictive distribution be continuous. This is
clearly not true in the case of the multinomial, as $M_{T_{n}}$ is always
discrete. It was proved in Evans and Jang (2011b), however, that when $k=1$ a
continuized version of $M_{T_{n}}$ can be constructed that yields the
consistency result in the binomial case. It is part of our purpose here to
establish the general consistency result for the multinomial. This turns out
to be much more difficult than the binomial case. As such the Appendix
contains the proof of the following result where $\Theta_{k,c},\Theta_{k,d}$
denote the sets of continuity and discontinuity points of $\pi.$

\begin{theorem}
\label{prior-in-conflict}Let $\pi$ be a prior on the probabilities $\theta
\in\Theta_{k}$ that satisfies the following conditions.\newline(A1) The prior
density is bounded above, that is, there exists $B>0$ such that $\pi
(\theta)\leq B$ for all $\theta\in\Theta_{k}$.\newline(A2) The prior density
is continuous almost surely with respect to volume measure, that is,
$\text{vol}(\Theta_{k,d})=0$.\newline(A3) The prior probability of each level
set of prior density is a null set with respect to the prior, that is,
$\Pi(\{\theta:\pi(\theta)=l\})=0$ for any $l\geq0$.\newline Then
(\ref{pvalue}) converges to $\Pi(\pi(\theta)\leq\pi(\theta^{\text{true}}))$
whenever $\theta^{\text{true}}\in\Theta_{k,c}$ as $n\rightarrow\infty.$
\end{theorem}

\noindent It is to be noted that, because of the discreteness, the value of
(\ref{pvalue}) is invariant under 1-1 transformations of $T$ and also under
reparametrizations. The theorem states that (\ref{pvalue}) converges to
$\Pi(\pi(\theta)\leq\pi(\theta^{\text{true}}))$ where $\pi$ is the prior on
the probabilities $(\theta_{1},\ldots,\theta_{k}).$ So if we reparameterized
and used some other prior to compute $m_{T},$ this has no effect on the limit.

Note that the theorem requires that $\Pi$ be a continuous probability measure
and also the prior density cannot be constant on a subregion of $\Theta_{k}$
having positive volume. For example, the theorem does not cover the uniform
prior on $\Theta_{k}$ although the result still holds there as (\ref{pvalue})
equals 1 as does $\Pi(\pi(\theta)\leq\pi(\theta^{\text{true}}))$ and there is
no need to check this prior$.$The following result also follows from the proof
of the theorem.

\begin{corollary}
\label{corprior-in-conflict}If $\pi$ satisfies (A1) and (A2) and is continuous
at $\theta^{\text{true}},$ then
\begin{align*}
&  \Pi(\pi(\theta)<\pi(\theta^{\text{true}}))\leq\liminf_{n\rightarrow\infty
}M_{T_{n}}(m_{T_{n}}(\mathbf{t})\leq m_{T_{n}}(T_{n}(x)))\\
&  \leq\limsup_{n\rightarrow\infty}M_{T_{n}}(m_{T_{n}}(\mathbf{t})\leq
m_{T_{n}}(T_{n}(x)))\leq\Pi(\pi(\theta)\leq\pi(\theta^{\text{true}})).
\end{align*}

\end{corollary}

\noindent So, if there is no prior-data conflict in the sense that $\Pi
(\pi(\theta)<\pi(\theta^{\text{true}}))$ is not small, then $M_{T_{n}%
}(m_{T_{n}}(\mathbf{t})\leq m_{T_{n}}(T_{n}(x)))$ for large $n$ will reflect this.

In the following examples we make use of a parameterization of the Dirichlet
referred to here as the \textit{concentration parameterization}. For the
Dirichlet$(\alpha_{1},\ldots,$\newline$\alpha_{k+1})$ distribution with all
$\alpha_{i}>1$ the mode is at $(\xi_{1},\ldots,\xi_{k+1})$ where $\xi
_{i}=(\alpha_{i}-1)/\tau$ and the \textit{concentration parameter}
$\tau=\alpha_{1}+\cdots+\alpha_{k+1}-(k+1).$ As $\alpha_{i}=1+\tau\xi_{i},$ it
is seen that the set of all Dirichlets with this mode is indexed by $\tau>0.$
The mean and variance of the $i$-th coordinate equal $(1+\tau\xi_{i}%
)/(\tau+k+1)$ and $(1+\tau\xi_{i})(\tau+k-\tau\xi_{i}))/(\tau+k+1)^{2}%
(\tau+k+2)$ which converge respectively to $\xi_{i}$ and 0 as $\tau
\rightarrow\infty,$ and so the distribution concentrates at the mode, and as
$\tau\rightarrow0$ the distribution converges to the uniform on the simplex.

Consider now the examples of .Section 2.\smallskip

\noindent\textbf{Example 5. }\textit{Checking the prior for the trine. }

For a single qubit, an experimenter without\ any prior knowledge could assign
a prior to the qubit state space that is uniform under the Hilbert-Schmidt
measure. When a trine measurement is performed on the qubit, this results in
the prior given by $\pi(\theta)\varpropto(1-\left(  \theta-c\right)
^{t}C\left(  \theta-c\right)  )^{1/2}$ when $\theta\in\Theta$ and $0$
otherwise, where $c$ and $C$ are as in Example 3. The change of variable
$\theta\rightarrow(r,\omega),$ where $\theta=c+C^{1/2}r^{1/2}(\cos\omega
,\sin\omega)^{t}$, has Jacobian proportional to $r^{1/2}.$ Therefore,
$\omega\sim U(0,2\pi)$ independent of $r\sim$ beta$(3/2,3/2)$ and this
provides an algorithm for generating from $\pi.$ In circumstances where $n$ is
modest, generating from $\pi$ and averaging the likelihood can be used to
compute the values $m_{T_{n}}(t_{1},t_{2},t_{3})$ needed for (\ref{pvalue}).
In this case $n$ is large so this is too inefficient due to the concentration
of the likelihood near the MLE over $\Theta_{k}.$ The posterior under the
uniform prior on $\Theta_{k}$ also concentrates near the MLE and so importance
sampling based on sampling from the Dirichlet$(t_{1}+1,t_{2}+1,t_{3}+1)$ is
used to estimate $m_{T_{n}}(t_{1},t_{2},t_{3})$ for each $(t_{1},t_{2},t_{3})$
in a sample drawn from $m_{T_{n}}.$ Sampling from $m_{T_{n}}$ is carried out
by generating $(\theta_{1},\theta_{2},\theta_{3})$ from $\pi$ and then
generating $(t_{1},t_{2},t_{3})\sim$ multinomial$(n,\theta_{1},\theta
_{2},\theta_{3}).$ The values of $m_{T_{n}}(t_{1},t_{2},t_{3})$ are compared
to $m_{T_{n}}(n_{1},n_{2},n_{3})$ to estimate (\ref{pvalue}).

This procedure was carried out for the entries in Table 1 with $10^{3}$ values
of $(t_{1},t_{2},t_{3})$ generated from $m_{T_{n}}$ and each value of
$m_{T_{n}}(t_{1},t_{2},t_{3})$ estimated using a sample of $10^{4}$ from the
relevant posterior based on $(t_{1},t_{2},t_{3})$. Prior-data conflict in this
example corresponds to the true value of $\theta$ lying near the boundary of
the respective set $\Theta.$ For the symmetric case (\ref{pvalue}) was
estimated as $0.87$ and for the asymmetric case (\ref{pvalue}) was estimated
as $0.15.$ These results were quite stable over different choices of
simulation sample sizes. So in neither case is there any indication of
prior-data conflict.

As the sample size are large, the importance samplers are quite concentrated.
For example, when estimating $m_{T_{n}}(n_{1},n_{2},n_{3}),$ the standard
deviations of the posterior distributions of the $\theta_{i},$ based on the
uniform prior on $\Theta_{2},$ are $(5.94\times10^{-3},5.28\times
10^{-3},5.13\times10^{-3})$ in the symmetric case. So to investigate the
sensitivity of the results, more diffuse importance samplers were considered
and this was achieved by taking lower values of $\tau$ in the concentration
parameterization of the Dirichlet. Here $\tau=n$ corresponds to using the
posterior based on the uniform prior on $\Theta_{2}$ as the importance sampler
and $\tau=0$ corresponds to using the uniform distribution on $\Theta_{2}$ as
the importance sampler which has standard deviations $(235.70\times
10^{-3},235.70\times10^{-3},235.70\times10^{-3}).$ Of course, as $\tau$ drops
it is to be expected that the efficiency of the importance sampling will
decrease. Even when $\tau=n/100,$ which has standard deviations $(57.76\times
10^{-3},51.51\times10^{-3},50.11\times10^{-3}),$ similar results were obtained
with the second decimal place in the estimate of (\ref{pvalue}) changing.
\smallskip

\noindent\textbf{Example 6. }\textit{Eliciting and checking the prior for
ordered probabilities. }

It is necessary to first provide an elicitation algorithm for a prior on
$\theta$ satisfying (\ref{ordered}). For this Lemma 1 helps considerably since
$(\theta_{1},\ldots,\theta_{k+1})^{t}=A_{k}(\omega_{1},\ldots,\omega
_{k+1})^{t}$ and so any prior on $\omega$ will induce a prior on $\theta$
satisfying (\ref{ordered}). Perhaps it is natural to choose a
Dirichlet$(\alpha_{1},\ldots,\alpha_{k+1})$ prior on $\omega$ but how should
the $\alpha_{i}$ be chosen? This of course depends upon what is known about
the $\theta_{i}$ and various elicitation algorithms can be considered.

Perhaps a natural approach is for the investigator to specify ordered
probabilities $\theta_{1}^{\ast}\geq\theta_{2}^{\ast}\geq\cdots\geq
\theta_{k+1}^{\ast}$ and then use a prior with mode at this point. By Lemma 1
this can be accomplished by a Dirichlet distribution with mode at $(\xi
_{1},\ldots,\xi_{k+1})^{t}=A_{k}^{-1}(\theta_{1}^{\ast},\ldots,\theta
_{k+1}^{\ast})^{t}$ and then, using the concentration parameterization, $\tau$
can be chosen to reflect belief in this mode. This requires, however, that the
$\theta_{i}^{\ast}$ satisfy (\ref{ordered}) as well as being a probability
distribution. The following result characterizes the choices of $\theta
_{i}^{\ast}$ that are equispaced as this seems like a somewhat natural choice
although there are many other possibilities that can be characterized in
similar ways.\smallskip

\noindent\textbf{Lemma 2}. The probabilities $\theta_{i}^{\ast}$ satisfying
(\ref{ordered}) are equispaced with $\theta_{i}^{\ast}=\theta_{1}^{\ast
}-(i-1)\delta$ iff $\theta_{1}^{\ast}=k\delta/2+1/(k+1)$ and $0\leq\delta
\leq2/k(k+1)$ and in this case $\xi_{i}=i\delta$ for $i=1,\ldots,k$ and
$\xi_{k+1}=1-k(k+1)\delta/2.$

\noindent\textbf{Proof}: By Lemma 1 the $\theta_{i}^{\ast}=\theta_{1}^{\ast
}-(i-1)\delta$ give probabilities satisfying (\ref{ordered}) for some
$\delta\geq0$ iff $\xi_{i}=i\delta$ for $i=1,\ldots,k$ and $\xi_{k+1}=\left(
k+1\right)  (\theta_{1}^{\ast}-k\delta)$ and $\theta_{1}^{\ast}\geq k\delta.$
Since $1=\sum_{j=1}^{k+1}\xi_{j}=k(k+1)\delta/2+\left(  k+1\right)
(\theta_{1}^{\ast}-k\delta)=\left(  k+1\right)  (\theta_{1}^{\ast}-k\delta/2)$
iff $\theta_{1}^{\ast}=k\delta/2+1/(k+1)$ and this satisfies $\theta_{1}%
^{\ast}\geq k\delta$ iff $\delta\leq2/k\left(  k+1\right)  .$\smallskip

\noindent Lemma 2 implies that $(\omega_{1},\ldots,\omega_{k+1})\sim
$Dirichlet$(1+\tau\delta,\ldots,1+k\tau\delta,1+\tau(1-k(k+1)\delta/2)).$ When
$\delta=0,$ then $\theta_{i}^{\ast}=1/(k+1)$ for all $i$ implying $\xi_{i}=0$
for $i=1,\ldots,k$ and $\xi_{k+1}=1$ so $(\omega_{1},\ldots,\omega_{k+1})\sim$
Dirichlet$(1,\ldots,1,1+\tau).$ When $\delta=2/k\left(  k+1\right)  ,$ then
$\theta_{i}^{\ast}=[2/(k+1)][1-(i-1)/k]$ for $i=1,\ldots,k$ implying $\xi
_{i}=2i/k(k+1)$ for $i=1,\ldots,k$ and $\xi_{k+1}=0$ so $(\omega_{1}%
,\ldots,\omega_{k+1})\sim$ Dirichlet$(1+2\tau/k(k+1),1+4\tau/k(k+1),\ldots
,1+2\tau/(k+1),1).$

It remains to determine $\tau.$ There are a number of ways to proceed but here
it is supposed that there is an interval $(l,u)$ such that it is believed all
the true probabilities lie in $(l,u)$ with virtual certainty, namely, the
prior satisfies%
\begin{equation}
\gamma\leq\Pi(l<\theta_{k+1},\theta_{1}<u)=\Pi\left(  l<\omega_{k+1}%
/(k+1),\sum_{i=1}^{k+1}\omega_{i}/i<u\right)  \label{virtcert}%
\end{equation}
where $\gamma$ is a large probability like $0.99.$ Since $l<\theta_{k+1}%
^{\ast},\theta_{1}^{\ast}<u,$ the right-hand side of (\ref{virtcert}) goes to
1 as $\tau\rightarrow\infty.$ Therefore, $\tau$ satisfying (\ref{virtcert}) is
easily found by simulation and the smallest such value of $\tau$ is preferable
as this implies the least concentration for the prior.

For the data of Example 4 and with $\delta=0$ then $\theta_{1}^{\ast}%
=\theta_{k+1}^{\ast}=1/18.$ Then $l=1/450,u=1/2$ are possible and the value
$\tau=2.85$ is the estimated smallest value satisfying (\ref{virtcert}) with
$\gamma=0.99.$ When $\delta=2/k\left(  k+1)\right)  ,$ then $\theta_{1}^{\ast
}=1/9$ and $\theta_{k+1}^{\ast}=0$ so $l=0,u=1/4$ is possible and the value
$\tau=16.5$ is the estimated smallest value satisfying (\ref{virtcert}) with
$\gamma=0.99.$ Now consider checking the prior corresponding to $\delta=0$ and
$l=1/450,u=1/2.$ Figure 2 is a plot of density histograms for the first four
probabilities based on a sample of $10^{5}$ from the full prior.%
\begin{figure}
[t]
\begin{center}
\includegraphics[
height=3.5155in,
width=3.5622in
]%
{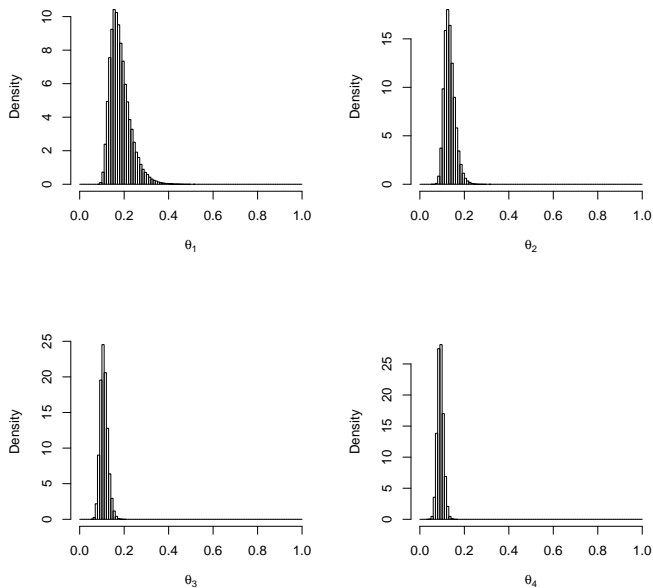}%
\caption{The marginal priors on $\theta_{1},\theta_{2},\theta_{3}$ and
$\theta_{4}$ when $\delta=0,l=1/450,u=1/2$ and $\tau=2.85$ in Example 6.}%
\end{center}
\end{figure}
Our approach to computing (\ref{pvalue}) is based on importance sampling. For
this particular data set $(f_{1}/n,\ldots,f_{18}/n)\in\Theta$ and an
importance sampler on $\Theta$ with this mode and values of $\tau\approx60$
produces reasonably stable estimates of $m(f_{1},\ldots,f_{18}).$ Note that
Lemma 1 plays a key role in obtaining such an importance sampler.

A problem arises, however, when computing values of $m(t_{1},\ldots,t_{18})$
necessary for estimating (\ref{pvalue}). Values of $(t_{1},\ldots,t_{18})$ are
obtained by generating $(\theta_{1},\ldots,\theta_{18})\sim\pi$ and then
$(t_{1},\ldots,t_{18})\sim$ multinomial$(n,\theta_{1},\ldots,\theta_{18}).$
When $n$ is small relative to dimension, as is the case here, then typically
$(t_{1}/n,\ldots,t_{18}/n)\notin\Theta$ and so the choice of importance
sampler is unclear. For example, an importance sampler such as the
Dirichlet$(t_{1}+1,\ldots,t_{18}+1),$ which has its mode at $(t_{1}%
/n,\ldots,t_{18}/n),$ virtually never generates points in $\Theta$ and so is
useless. A better approach is to use an importance sampler based on Lemma 1
where the Dirichlet on $(\omega_{1},\ldots,\omega_{k+1})$ has its mode at
$A^{-1}(\theta_{1}^{\ast\ast},\ldots,\theta_{18}^{\ast\ast})^{t}$ where
$(\theta_{1}^{\ast\ast},\ldots,\theta_{18}^{\ast\ast})$ is the convex
combination of the prior mode and $(t_{1}/n,\ldots,t_{18}/n)$ that just
satisfies being in $\partial\Theta.$ This always generates points inside
$\Theta$ and should at least somewhat mimic the integrand over $\Theta$
provided the concentration $\tau$ is not chosen too large or too small. Here a
representative $(t_{1},\ldots,t_{18})$ was selected and $\tau$\ chosen such
that both smaller and larger values lead to smaller estimates of
(\ref{pvalue}). When this was carried out on this example the value of
(\ref{pvalue}) was estimated as $0.36$ which indicates there is no prior-data
conflict. This makes sense as the naive estimate $(f_{1}/n,\ldots,f_{18}/n)$
is not only in $\Theta$ but also satisfies the bounds. Suppose instead the
data $f=(35,29,20,145,96,11,4,4,4,3,3,2,2,1,1,1,1,1)$ was observed that
clearly violates the monotonicity. In this case (\ref{pvalue}) equals 0 and so
prior-data conflict was detected as is correct. Also, with the original data
and the prior determined by $l=0,u=0.2$ satisfying (\ref{virtcert}), then
(\ref{pvalue}) equals 0 again indicating a definite prior-data conflict.

One way to avoid the computational problems encountered estimating
(\ref{pvalue}) is to collect more data as the distribution of $(t_{1}%
/n,\ldots,t_{k+1}/n)$ becomes degenerate at a point in $\Theta$ as
$n\rightarrow\infty$ and so will be in $\Theta$ for all $n$ large enough.
Another approach to avoiding these problems is to reduce dimension by grouping
as was done in Example 4. Intuitively, as the ratio of dimension to sample
size decreases, the values of the generated relative frequencies are more
likely to be in the relevant parameter space and this was confirmed by a
simulation experiment. Implementing the importance sampling for estimating the
prior predictive densities of the reduced problem, however, requires the
computation of the marginal prior density as this is not in closed form. Since
this would be required at each generated value from the importance sampler,
the computational advantage is largely negated. The prior density of the full
parameter $\theta$ is easily obtained via (\ref{basiceq}) so this is not an
issue in that case.

As such, an alternative approach is considered based on the original
elicitation algorithm and which could be called \textit{marginalizing the
elicitation}. Rather than trying to marginalize the full prior, consider being
presented with the reduced problem and then applying the elicitation algorithm
to that problem. This will not result in a prior that is the marginal of the
full prior but surely checking this prior for conflict with the data is also
assessing whether or not the information being used to choose the prior is
appropriate. For example, the original elicitation led to the inequality
$u\geq\theta_{1}\geq\cdots\geq\theta_{k+1}\geq l$ holding with virtual
certainty and recall that necessarily $\theta_{i}\leq1/i.$ So, if cells are
grouped in pairs to maintain the monotonicity and to make the best use of the
bounds, then supposing $k+1$ is even, $u+1/(1+(k+1)/2)\geq\theta_{1}%
+\theta_{1+(k+1)/2}\geq\theta_{2}+\theta_{2+(k+1)/2}\geq\cdots\geq
\theta_{(k+1)/2}+\theta_{k+1}\geq l+1/(k+1).$ If $k+1$ is odd then the last
group can consist of $\theta_{k+1}$ by itself and the lower bound doesn't
change. So for even modest $k$ the bounds will not increase by much and
clearly this idea can be extended to groups of 3, 4 etc. Lemma 1 can then be
used, as in the full problem, to obtain the prior for the parameters for the
grouped problem. Supposing there are $m$ groups and $(t_{1}^{\text{red}%
},\ldots,t_{m}^{\text{red}})$ denotes a value generated from the prior
predictive for the reduced problem, our recommendation is that this reduction
be continued until a reasonable proportion of the values $(t_{1}^{\text{red}%
}/n,\ldots,t_{m}^{\text{red}}/n)$ lie in the reduced parameter space. When
this is the case even those that lie outside should be close to the relevant
set which will improve the quality of the importance sampling. Note that this
does not imply that the observed data when grouped has to lie in the reduced
parameter space, as there may indeed be prior-data conflict, but because the
model has been accepted, it seems likely that this point will be either in or
close to this set.

For the model considered here, with $10^{4}$ values of $(t_{1}^{\text{red}%
}/n,\ldots,t_{m}^{\text{red}}/n)$ generated from the prior, the following
values of $(m,p_{m})$ were obtained where $p_{m}$ is the proportion that lay
inside the relevant parameter space: $(18,0.00),(9,0.04),$\newline%
$(6,0.25),(5,0.44),(4,0.61),(3,0.78)$ and $(2,0.93).$ The values $0.42$ and
$0.44$ were obtained for (\ref{pvalue}) when $m=9$ and $m=6,$ respectively. So
one can feel fairly confident that the elicited information is not in conflict
with the data.

If the model and prior are deemed acceptable, then computations based on the
posterior are required for inference and these can proceed using importance
sampling as described. Also, a Gibbs sampling approach is available. Putting
$\theta_{\neq i}=(\theta_{1},\ldots,\theta_{i-1},\theta_{i+1},\ldots
,\theta_{k})$ with $\alpha_{0}=1,\theta_{0}=1$ when $i=1,\ldots,k-1,$ then
$\theta_{i}\,|\,\theta_{\neq i}$ has density proportional to $\theta
_{i}^{f_{i}}(1-\theta_{i}-\theta_{\neq i})^{f_{k+1}+\alpha_{k+1}-1}%
(\theta_{i-1}-\theta_{i})^{\alpha_{i-1}-1}(\theta_{i}-\theta_{i+1}%
)^{\alpha_{i}-1}(\theta_{i}+\theta_{k}-1+\theta_{\neq i})^{\alpha_{k}-1}$ for
$\theta_{i}$ in the fixed interval $[\max\{\theta_{i+1},1-\theta_{\neq
i}-\theta_{k}\},\min\{\theta_{i-1},1-\theta_{\neq i}\}]$ and for $i=k,$ then
$\theta_{k}\,|\,\theta_{\neq k}$ has density proportional to $\theta
_{k}^{f_{k}}(1-\theta_{k}-\theta_{\neq k})^{f_{k+1}+\alpha_{k+1}-1}%
(\theta_{k-1}-\theta_{k})^{\alpha_{k-1}-1}(2\theta_{k}-1+\theta_{\neq
k})^{\alpha_{k}-1}$ for $\theta_{k}$ in the fixed interval $[(1-\theta_{\neq
k})/2,\min\{\theta_{k-1},1-\theta_{\neq k}\}].$ Sampling from these densities
can proceed by generating $\theta_{i}/(1-\theta_{\neq i})\,|\,\theta_{\neq i}$
as a truncated beta$(f_{i}+1,f_{k+1}+\alpha_{k+1}),$ which accounts for the
first two factors, and then using the envelope methods described in Evans and
Swartz (1998) to account for the remaining factors.

The ZM$_{k}$ model was also considered but a significant problem with this
family remains unresolved. In particular, it is unclear how to elicit a prior
on $\left(  \alpha,\beta\right)  $ and this is because the interpretation of
these parameters is not obvious. Also, the ZM$_{k}$\ family imposes some
sharper constraints on the probabilities than hold generally. For example, the
maximum probabilities over all $\left(  \alpha,\beta\right)  $ for $i=1,2,3,4$
are $1.00,0.33,0.26,0.25,$ respectively, which contrasts with
$1.00,0.50,0.33,0.25$ for the general model for ordered probabilities. So,
given that it is much easier to use and interpret, the general model for
ordered probabilities is recommended over the ZM$_{k}$ model.

\section{Conclusions}

Constrained multinomial models arise in a number of interesting contexts and
pose some unique challenges. The emphasis here has been on checking these
models and checking for prior-data conflict. Issues associated with inference
will be considered elsewhere.\ A significant consistency theorem has been
established for the check on the prior. As a particular application a general
model for ordered probabilities has been developed, together with an
elicitation algorithm for a prior, and the results of the paper applied. Also,
one of a variety of constrained multinomial models used in quantum state
estimation has been used to illustrate the methodology.

\section*{Acknowledgements}

This work is funded by the Singapore Ministry of Education (partly through the
Academic Research Fund Tier 3 MOE2012-T3-1-009) and the National Research
Foundation of Singapore.

\section*{Appendix}

The proof of the theorem proceeds via a number of lemmas which are
individually of some interest. As with the binomial case it is necessary to
construct a continuous probability measure that is essentially equivalent to
$M_{T_{n}}$.

Let $D_{k,n}=\{\mathbf{n}=(n_{1},\ldots,n_{k}):n_{i}\in%
\mathbb{N}
_{0},n_{1}+\cdots+n_{k}\leq n\}$ denote the set of possible values of $T_{n}.$
Now construct a set of disjoint sets that cover $\Theta_{k}$ and are indexed
by $\mathbf{n\in}D_{k,n}$ such that $\mathbf{n/}n$ is in the set that
$\mathbf{n}$ indexes. Let $\Theta_{k}(\mathbf{n})=%
{\textstyle\prod_{i=1}^{k}}
[n_{i}/n-1/2n,n_{i}/n+1/2n)$ and the $\Theta_{k}(\mathbf{n})$ are disjoint,
$\mathbf{n/}n$ is the center of $\Theta_{k}(\mathbf{n}),$ vol$(\Theta
_{k}(\mathbf{n}))=n^{-k}$ and $\Theta_{k}\subset\cup_{\mathbf{n\in}D_{k,n}%
}\Theta_{k}(\mathbf{n})\downarrow\bar{\Theta}_{k}$ as $n\rightarrow\infty.$

For $\mathbf{r\in}\cup_{\mathbf{n\in}D_{k,n}}\Theta_{k}(\mathbf{n})$ define
\[
m_{n}^{\ast}(\mathbf{r})=n^{k}m_{T_{n}}(\mathbf{n}(\mathbf{r}))
\]
where $\mathbf{n}(\mathbf{r})\in D_{k,n}$ is such that $\mathbf{r\in}%
\Theta_{k}(\mathbf{n}(\mathbf{r}))$. Note that for $\mathbf{r\in}\Theta
_{k}(\mathbf{n}(\mathbf{r})),$ then $\mathbf{n}(\mathbf{r})=(n_{1}%
(\mathbf{r)},\ldots,n_{k}(\mathbf{r)})$ where $n_{i}(\mathbf{r)}=\lfloor
nr_{i}+1/2\rfloor$ and $\lfloor x\rfloor$ is the biggest integer that is not
greater than $x.$ Also, define $c_{n}(\mathbf{r)=n}(\mathbf{r})/n$ which is
the center of the $k$-cell containing $\mathbf{r.}$ The following result then holds.

\begin{lemma}
\label{oblemm}(i) $||\mathbf{r-}c_{n}\mathbf{(\mathbf{r})}||\leq\sqrt{k}/n$ so
$c_{n}\mathbf{(\mathbf{r})}\rightarrow\mathbf{r}$ as $n\rightarrow\infty.$
(ii) $m_{n}^{\ast}$ is constant on $\Theta_{k}(\mathbf{n})$ and takes the
value $n^{k}m_{T_{n}}(\mathbf{n})$ there. (iii) $m_{n}^{\ast}$ is the density
of an absolutely continuous probability measure $M_{n}^{\ast}$ where
$M_{n}^{\ast}(\Theta_{k}(\mathbf{n}))=m_{T_{n}}(\mathbf{n})$ and
\[
M_{T_{n}}(\{\mathbf{t}:m_{T_{n}}(\mathbf{t})\leq m_{T_{n}}(T_{n}%
(x))\})=M_{n}^{\ast}(\{\mathbf{r}:m_{n}^{\ast}(\mathbf{r})\leq m_{n}^{\ast
}(T_{n}(x)/n)\}).
\]

\end{lemma}

\noindent\textbf{Proof}: Parts (i) and (ii) are obvious. For (iii) we have
$M_{T_{n}}(\{\mathbf{t}:m_{T_{n}}(\mathbf{t})\leq m_{T_{n}}(T_{n}%
(x))\})=M_{T_{n}}(\{n(\mathbf{r}):m_{n}^{\ast}(\mathbf{r})\leq m_{n}^{\ast
}(T_{n}(x)/n)\}).\smallskip$

\noindent So indeed $M_{n}^{\ast}$ is a continuized version of $M_{n}.$

For $s>0,$ define%
\begin{align*}
G_{n,s}(\mathbf{r})  &  =\{\mathbf{\theta}\in\Theta_{k}:\sum_{i=1}^{k+1}%
c_{ni}\mathbf{(\mathbf{r})}\log(\theta_{i}/c_{ni}\mathbf{(\mathbf{r}%
)})>-[(k+1+s)\log(n)]/n\}\\
&  =\{\mathbf{\theta}\in\Theta_{k}:\text{KL}(c_{n}(\mathbf{r)\,||\,\theta
)<}[(k+1+s)\log(n)]/n\},
\end{align*}
a Kullback-Leibler neighborhood of $c_{n}(\mathbf{r).}$ We have the following
result where $B_{\epsilon}(\mathbf{r})$ is the ball of radius $\epsilon$
centered at $\mathbf{r}.$

\begin{lemma}
\label{onlyrlemm}If $\mathbf{r\in}\Theta_{k},$ then, for $\epsilon>0,$ there
exists $N$ such that for all $n>N,\mathbf{r\in\,}G_{n,s}(\mathbf{r})\subset
B_{\epsilon}(\mathbf{r}),$ so $G_{n,s}(\mathbf{r})\longrightarrow
\{\mathbf{r\}}$ as $n\rightarrow\infty.$
\end{lemma}

\noindent\textbf{Proof}: Note that $r_{i}\neq0,$ so for $n$ large enough
$\lfloor nr_{i}+1/2\rfloor\geq1\ $and this holds for all $i$. Since
$nr_{i}+1/2\geq$ $\lfloor nr_{i}+1/2\rfloor$ and $\log(nr_{i}/(nr_{i}%
+1/2))<0$,%
\begin{align*}
&  \sum_{i=1}^{k+1}c_{n,i}(\mathbf{r)}\log\left(  \frac{r_{i}}{c_{n,i}%
}\right)  =\sum_{j=1}^{k+1}\left(  \frac{\lfloor nr_{i}+1/2\rfloor}{n}\right)
\log\left(  \frac{nr_{i}}{\lfloor nr_{i}+1/2\rfloor}\right)  \\
&  \geq\sum_{i=1}^{k+1}\left(  \frac{\lfloor nr_{i}+1/2\rfloor}{n}\right)
\log\left(  \frac{nr_{i}}{nr_{i}+1/2}\right)  \geq\sum_{i=1}^{k+1}\left(
\frac{nr_{i}+1/2}{n}\right)  \log\left(  \frac{nr_{i}}{nr_{i}+1/2}\right)  \\
&  =\sum_{i=1}^{k+1}\left(  \frac{nr_{i}+1/2}{n}\right)  \log\left(
1-\frac{1}{2(nr_{i}+1/2)}\right)  \\
&  =-\frac{1}{2n}\sum_{i=1}^{k+1}\left\{  \sum_{j=1}^{\infty}\frac{\left[
2(nr_{i}+1/2)\right]  ^{-j+1}}{j}\right\}  \geq-\frac{1}{2n}\sum_{i=1}%
^{k+1}\left\{  \sum_{j=1}^{\infty}\left[  2(nr_{i}+1/2)\right]  ^{-j+1}%
\right\}  \\
&  =-\frac{1}{2n}\sum_{i=1}^{k+1}\frac{1}{1-1/2(nr_{i}+1/2)}=-\frac{1}{2n}%
\sum_{i=1}^{k+1}\left(  1+\frac{1}{2nr_{i}}\right)  =-\frac{k+1}{2n}-\frac
{1}{4n^{2}}\sum_{i=1}^{k+1}\frac{1}{r_{i}}%
\end{align*}
and clearly there is an $N$ such that this is bounded below by $-[(k+1+s)\log
(n)]/n$ for all $n>N.$ Therefore, $\mathbf{r\in}G_{n,s}(\mathbf{r})$ for all
$n$ large enough. If $||\mathbf{r-\theta}||=\epsilon>0,$ then by Lemma
\ref{oblemm} (i), $||c_{n}\mathbf{(\mathbf{r})-\theta}||>\epsilon/2$ for all
$n$ large enough. Since $\sum_{j=1}^{k+1}c_{nj}\mathbf{(\mathbf{r})}%
\log(\theta_{j}/c_{nj}\mathbf{(\mathbf{r})})$ is continuous in $c_{n}%
\mathbf{(\mathbf{r})}$, bounded above by 0 and 0 iff $c_{n}\mathbf{(\mathbf{r}%
)=\theta,}$ then
\[
\sum_{j=1}^{k+1}c_{nj}\mathbf{(\mathbf{r})}\log(\theta_{j}/c_{nj}%
\mathbf{(\mathbf{r})})\rightarrow\delta<0.
\]
From this we conclude that $\mathbf{\theta\notin}G_{n,s}(\mathbf{r})$ for all
$n$ large enough and this completes the proof.$\smallskip$

If $\mathbf{r}\mathbf{\notin}\bar{\Theta}_{k},$ then clearly $m_{n}^{\ast
}(\mathbf{r})\rightarrow\pi(\mathbf{r})=0.$ It is now proved that $m_{n}%
^{\ast}(\mathbf{r})\rightarrow\pi(\mathbf{r})$ when $\mathbf{r\in}\Theta_{c}.$
The following identity is useful, namely, when $\mathbf{n}\in D_{k,n}%
\mathbf{,}$ then%
\begin{equation}
\int_{\Theta_{k}}n!\prod_{j=1}^{k+1}\left(  \theta_{j}^{n_{j}}/n_{j}!\right)
\,d\mathbf{\theta}=n!/(n+k)!.\label{id1}%
\end{equation}

\begin{lemma}
\label{lem:mnstar-approx}If $\pi$ satisfies assumptions A1 and A2 then, for
$\mathbf{r\in}\Theta_{k}$ and any $s>0$,
\[
\left\vert m_{n}^{\ast}(\mathbf{r})-\frac{n!n^{k}}{(n+k)!}\pi(\mathbf{r}%
)\right\vert \leq\frac{1}{n^{s}}+\sup_{\mathbf{\theta}\in G_{n,s}%
(\mathbf{r})\cap\Theta_{k,c}}|\pi(\mathbf{\theta})-\pi(\mathbf{r})|
\]
when $n>(1+Bk!)^{2}.$
\end{lemma}

\noindent\textbf{Proof}: We abbreviate $G_{n,s}(\mathbf{r})$ to $G_{n,s}$ and
let $F_{n,s}=\Theta_{k}-G_{n,s}=\{\mathbf{\theta}:\sum_{j=1}^{k+1}c_{n,j}%
\log(\theta_{j}/c_{n,j})\leq-[(k+1+s)\log(n)]/n\}$. Then
\begin{align*}
&  \left\vert m_{n}^{\ast}(\mathbf{r})-\frac{n!n^{k}}{(n+k)!}\pi
(\mathbf{r})\right\vert =\left\vert n^{k}\int_{\Theta_{k}}n!\prod_{j=1}%
^{k+1}\frac{\theta_{j}^{n_{j}}}{n_{j}!}\Pi(d\mathbf{\theta})-n^{k}%
\pi(\mathbf{r})\int_{\Theta_{k}}n!\prod_{j=1}^{k+1}\frac{\theta_{j}^{n_{j}}%
}{n_{j}!}\,d\mathbf{\theta}\right\vert \\
&  \leq n^{k}\int_{F_{n,s}}n!\prod_{j=1}^{k+1}\frac{\theta_{j}^{n_{j}}}%
{n_{j}!}\,\Pi(d\mathbf{\theta})+\pi(\mathbf{r})n^{k}\int_{F_{n,s}}%
n!\prod_{j=1}^{k+1}\frac{\theta_{j}^{n_{j}}}{n_{j}!}\,d\mathbf{\theta}+\\
&  \hspace{1in}n^{k}\int_{G_{n,s}}n!\prod_{j=1}^{k+1}\frac{\theta_{j}^{n_{j}}%
}{n_{j}!}|\pi(\mathbf{\theta})-\pi(\mathbf{r})|\,d\mathbf{\theta}%
=I_{n,1}+I_{n,2}+I_{n,3}.
\end{align*}

\noindent We will find upper bounds for the three terms $I_{n,1}%
,I_{n,2},I_{n,3}$.

First we show that $I_{n,1}\leq n^{-s-1/2}$ for all $n$. For $\mathbf{\theta
}\in F_{n,s}$, putting $n_{j}=nc_{n,j}$,%
\[
\sum_{j=1}^{k+1}n_{j}\log(\theta_{j})\leq\sum_{j=1}^{k+1}n_{j}\log
(n_{j})-(n+k+1+s)\log(n),
\]
and the probability function satisfies
\[
n!\prod_{j=1}^{k+1}\frac{1}{n_{j}!}\times\prod_{j=1}^{k+1}\theta_{j}^{n_{j}%
}\leq n!\prod_{j=1}^{k+1}\frac{1}{n_{j}!}\times\frac{1}{n^{n}}\prod
_{j=1}^{k+1}n_{j}^{n_{j}}\times\frac{1}{n^{k+1+s}}=\frac{n!}{n^{n}}%
\prod_{j:n_{j}>0}\frac{n_{j}^{n_{j}}}{n_{j}!}\times\frac{1}{n^{k+1+s}}.
\]
When $\max(n_{1},\ldots,n_{k+1})=n$, that is, $n_{j}=n$ for some $j$ and
$n_{i}=0$ for all $i\not =j$,
\[
\frac{n!}{n^{n}}\prod_{j:n_{j}>0}\frac{n_{j}^{n_{j}}}{n_{j}!}=\frac{n!}{n^{n}%
}\frac{n^{n}}{n!}=1\leq n^{1/2}.
\]
Assume $\max(n_{1},\ldots,n_{k+1})<n$. Then by the Robbins (1955) result on
Stirling's approximation,
\[
1/(12n+1)<\log(n!)-\frac{1}{2}\log(2\pi n)+n\log(n)-n<1/12n
\]
for all $n>0$, and so
\begin{align*}
\frac{n!}{n^{n}}\prod_{j:n_{j}>0}\frac{n_{j}^{n_{j}}}{n_{j}!} &  \leq(2\pi
n)^{1/2}e^{-n}e^{1/12n}\prod_{j:n_{j}>0}\frac{e^{n_{j}}e^{-1/(12n_{j}+1)}%
}{(2\pi n_{j})^{1/2}}\\
&  \leq\frac{n^{1/2}}{(2\pi)^{1/2}}\prod_{j:n_{j}>0}\frac{1}{n_{j}^{1/2}%
}<n^{1/2}%
\end{align*}
since $0<\exp\{1/12n-\sum_{j:n_{j}>0}1/(12n_{j}+1)\}<1.$ Then for any prior
$\Pi$
\begin{align*}
n^{k}\int_{F_{n,s}}n!\prod_{j:n_{j}>0}\frac{\theta_{j}^{n_{j}}}{n_{j}!}%
\,\Pi(d\mathbf{\theta}) &  \leq n^{k}\int_{F_{n,s}}n^{1/2}\frac{1}{n^{k+1+s}%
}\,\Pi(d\mathbf{\theta})\leq\frac{1}{n^{s+1/2}}\Pi(F_{n,s})\\
&  \leq\frac{1}{n^{s+1/2}}.
\end{align*}
Hence $I_{n,1}\leq n^{-(s+1/2)}$ regardless of prior. \medskip

The Dirichlet$(1,\ldots,1)$ has density equal to $k!$ on $\Theta.$ Applying
the above argument with $\Pi$ equal to this prior implies $I_{n,2}\leq
\pi(\mathbf{r})k!n^{-s-1/2}$. By A2%
\[
I_{n,3}=n^{k}\int_{G_{n,s}\cap\Theta_{k,c}}n!\prod_{j=1}^{k+1}\frac{\theta
_{j}^{n_{j}}}{n_{j}!}|\pi(\mathbf{\theta})-\pi(\mathbf{r})|\,d\mathbf{\theta}%
\]
which, using (\ref{id1}), implies
\[
I_{n,3}\leq\frac{n!n^{k}}{(n+k)!}\times\sup_{\mathbf{\theta}\in G_{n,s}%
\cap\Theta_{k,c}}|\pi(\mathbf{\theta})-\pi(\mathbf{r})|\leq\sup
_{\mathbf{\theta}\in G_{n,s}\cap\Theta_{k,c}}|\pi(\mathbf{\theta}%
)-\pi(\mathbf{r})|
\]
since $n!n^{k}/(n+k)!\leq1.$ Finally, for $n\geq(1+Bk!)^{2}$ and using A1,
\begin{align*}
&  I_{n,1}+I_{n,2}+I_{n,3}\\
&  \leq\frac{1+Bk!}{n^{1/2}}\frac{1}{n^{s}}+\sup_{\mathbf{\theta}\in
G_{n,s}\cap\Theta_{k,c}}|\pi(\mathbf{\theta})-\pi(\mathbf{r})|\leq\frac
{1}{n^{s}}+\sup_{\mathbf{\theta}\in G_{n,s}\cap\Theta_{k,c}}|\pi
(\mathbf{\theta})-\pi(\mathbf{r})|.
\end{align*}
and the lemma is proved.$\smallskip$

\begin{corollary}
\label{cor:mnstar-approx}$m_{n}^{\ast}\rightarrow\pi$ as $n\rightarrow\infty$
where the convergence is almost sure with respect to volume measure.
\end{corollary}

\noindent\textbf{Proof}: If $\mathbf{r\notin}\bar{\Theta}_{k},$ then
$m_{n}^{\ast}(\mathbf{r})\rightarrow\pi(\mathbf{r})=0.$ Now suppose
$\mathbf{r\in}\Theta_{k,c}.$ For $\epsilon>0$ there exists $\delta$ such that
$\mathbf{\theta\in}B_{\delta}(\mathbf{r})$ implies $|\pi(\mathbf{\theta}%
)-\pi(\mathbf{r})|<\epsilon/3.$ By Lemma \ref{onlyrlemm} there exists
$N>(1+Bk!)^{2}$ such that for all $n>N,$ then $G_{n,s}(\mathbf{r)\subset
}B_{\delta}(\mathbf{r}),n^{-s}<\epsilon/3$ and $1-n!n^{k}/(n+k)!<\epsilon/3B.$
Therefore, by Lemma \ref{lem:mnstar-approx}
\[
\left\vert m_{n}^{\ast}(\mathbf{r})-\pi(\mathbf{r})\right\vert \leq\left\vert
m_{n}^{\ast}(\mathbf{r})-\frac{n!n^{k}}{(n+k)!}\pi(\mathbf{r})\right\vert
+\left\vert \frac{n!n^{k}}{(n+k)!}\pi(\mathbf{r})-\pi(\mathbf{r})\right\vert
<\epsilon
\]
establishing convergence everywhere except on $\Theta_{k,d}\cup\partial
\Theta_{k}$ which has volume $0$. $\smallskip$

The following technical result is required for the proof of the theorem.

\begin{lemma}
\label{klnbdlemm}If $\pi$ is continuous at $\theta^{\text{true}}$, then as
$n\rightarrow\infty$%
\[
\sup_{\mathbf{\theta}\in G_{n,s}(T_{n}(x)/n)\cap\Theta_{k,c}}\left\vert
\pi(\mathbf{\theta})-\pi(\theta^{\text{true}})\right\vert \rightarrow0.
\]

\end{lemma}

\noindent\textbf{Proof}: Since $\pi$ is continuous at $\theta^{\text{true}},$
there exists $\delta>0$ such that, if $\mathbf{\theta}\in B_{\delta}%
(\theta^{\text{true}}),$ then $\left\vert \pi(\mathbf{\theta})-\pi
(\theta^{\text{true}})\right\vert <\epsilon.$ If $\mathbf{\theta}\in
G_{n,s}(T_{n}(x)/n),$ then the Kullback-Csiszar-Kemperman inequality (see
Devroye (1987), Theorem 1.4) says
\[
\sum_{j=1}^{k+1}|\theta_{j}-T_{n,j}/n|\leq\left[  2\text{KL}(T_{n}%
(x)/n\mathbf{\,||\,\theta})\right]  ^{1/2}<\left[  2n^{-1}(k+1+s)\log
n\right]  ^{1/2}%
\]
which implies that
\[
\sum_{j=1}^{k+1}|\theta_{j}-\theta_{j}^{\text{true}}|\leq(2n^{-1}(k+1+s)\log
n)^{1/2}+\sum_{j=1}^{k+1}|\theta_{j}^{\text{true}}-T_{n,j}/n|.
\]
Therefore, the almost sure convergence $T_{n}/n\rightarrow\theta^{\text{true}%
}$ implies that there exists $N$ such that for all $n>N,$ if $\mathbf{\theta
}\in G_{n,s}(T_{n}(x)/n),$ then $||\mathbf{\theta-}\theta^{\text{true}%
}||<\delta.$ This proves the lemma. $\smallskip$

The main result is now established.$\smallskip$

\noindent\textbf{Proof of Theorem }\ref{prior-in-conflict}: Fix $0<\eta<1$
small. Under A3, $\pi(\mathbf{\theta})$ has a continuous distribution when
$\mathbf{\theta\sim}\pi.$ Therefore, there exists $\epsilon>0$ such that
$\Pi(\{\mathbf{\theta}:|\pi(\mathbf{\theta})-\pi(\theta^{\text{true}}%
)|\leq\epsilon\})<\eta$.

Define $H_{n}=\{\mathbf{r}\in\Theta_{k,c}:\sup_{\mathbf{\theta}\in
G_{n,s}(\mathbf{r})\cap\Theta_{k,c}}|\pi(\mathbf{\theta})-\pi(\mathbf{r}%
)|<\epsilon/6\}.$ By Lemma \ref{onlyrlemm} the diameter of $G_{n,s}%
(\mathbf{r})$ shrinks to zero. If $\mathbf{r}\in\Theta_{k,c}$, then
$\sup_{\mathbf{\theta}\in G_{n,s}(\mathbf{r})\cap\Theta_{k,c}}|\pi
(\mathbf{\theta})-\pi(\mathbf{r})|\rightarrow0$ as $n\rightarrow\infty$, so
there exists $N(\mathbf{r},\epsilon)>0$ such that $\mathbf{r}\in H_{n}$ for
all $n\geq N(\mathbf{r},\epsilon).$ This implies $H_{n}\rightarrow\Theta
_{k,c}$ and so $\Pi(H_{n})\rightarrow\Pi(\Theta_{k,c})=1$ as $n\rightarrow
\infty$.

From the continuity of $\pi$ at $\theta^{\text{true}}$, there exists
$\delta>0$ such that $|\pi(\mathbf{\theta})-\pi(\theta^{\text{true}%
})|<\epsilon/12$ for any $\mathbf{\theta}\in B_{\delta}(\theta^{\text{true}}%
)$. The strong law of large numbers implies $T_{n}(x)/n\rightarrow
\theta^{\text{true}}$ almost surely so there exits $N_{1}$ such that
$||T_{n}(x)/n-\theta^{\text{true}}||<\delta$ for all $n>N_{1}$. Therefore, if
$n>N_{1},$ then $|\pi(T_{n}(x)/n)-\pi(\theta^{\text{true}})|<\epsilon/12.$
Also, by Lemma \ref{klnbdlemm} there exists $N_{2}$ such that for $n>N_{2},$
then
\[
\sup_{\mathbf{\theta}\in G_{n,s}(T_{n}(x)/n)\cap\Theta_{k,c}}\left\vert
\pi(\mathbf{\theta})-\pi(\theta^{\text{true}})\right\vert <\epsilon/6
\]
and $n^{-s}<\epsilon/6.$ Therefore, using Lemma \ref{lem:mnstar-approx}, for
all $n>N_{3}=\max\{(1+B_{m}k!)^{2},N_{1},N_{2}\},$%
\begin{align*}
&  \left\vert m_{n}^{\ast}\left(  \frac{T_{n}(x)}{n}\right)  -\frac{n!n^{k}%
}{(n+k)!}\pi(\theta^{\text{true}})\right\vert \\
&  \leq\left\vert m_{n}^{\ast}\left(  \frac{T_{n}(x)}{n}\right)
-\frac{n!n^{k}}{(n+k)!}\pi\left(  \frac{T_{n}(x)}{n}\right)  \right\vert
+\frac{n!n^{k}}{(n+k)!}\left\vert \pi\left(  \frac{T_{n}(x)}{n}\right)
-\pi(\theta^{\text{true}})\right\vert \\
&  \leq\frac{1}{n^{s}}+\sup_{\mathbf{\theta}\in G_{n,s}(T_{n}(x)/n)\cap
\Theta_{k,c}}\left\vert \pi(\mathbf{\theta})-\pi\left(  \frac{T_{n}(x)}%
{n}\right)  \right\vert +\\
&  \hspace{2in}\frac{n!n^{k}}{(n+k)!}\left\vert \pi\left(  \frac{T_{n}(x)}%
{n}\right)  -\pi(\theta^{\text{true}})\right\vert \\
&  \leq\frac{1}{n^{s}}+\sup_{\mathbf{\theta}\in G_{n,s}(T_{n}(x)/n)\cap
\Theta_{k,c}}\left\vert \pi(\mathbf{\theta})-\pi\left(  \mathbf{\theta}%
_{true}\right)  \right\vert +\\
&  \hspace{2in}\left(  1+\frac{n!n^{k}}{(n+k)!}\right)  \left\vert \pi
(\theta^{\text{true}})-\pi\left(  \frac{T_{n}(x)}{n}\right)  \right\vert \\
&  \leq\epsilon/2.
\end{align*}
Note that this implies that for all $n>N_{3},$%
\begin{equation}
m_{n}^{\ast}(T_{n}(x)/n)\leq\frac{n!n^{k}}{(n+k)!}\pi(\theta^{\text{true}%
})+\epsilon/2\leq\pi(\theta^{\text{true}})+\epsilon/2.\label{ineq1}%
\end{equation}

The prior-data conflict probability satisfies
\[
M_{n}^{\ast}(\{\mathbf{r}:m_{n}^{\ast}(\mathbf{r})\leq m_{n}^{\ast}%
(T_{n}(x)/n)\})=M_{n}^{\ast}(\{\mathbf{r}:m_{n}^{\ast}(\mathbf{r})\leq
m_{n}^{\ast}(T_{n}(x)/n)\}\cap\Theta_{k,c})
\]
since $M_{n}^{\ast}$ is absolutely continuous and $\Theta_{k,d}$ has volume
measure 0. Also,%
\begin{align*}
M_{n}^{\ast}(\{\mathbf{r}:m_{n}^{\ast}(\mathbf{r})  & \leq m_{n}^{\ast}%
(T_{n}(x)/n)\}\cap\Theta_{k,c})=\\
M_{n}^{\ast}(\{\mathbf{r}  & :m_{n}^{\ast}(\mathbf{r})\leq m_{n}^{\ast}%
(T_{n}(x)/n)\}\cap H_{n})+\\
& M_{n}^{\ast}(\{\mathbf{r}:m_{n}^{\ast}(\mathbf{r})\leq m_{n}^{\ast}%
(T_{n}(x)/n)\}\cap H_{n}^{c}\cap\Theta_{k,c})
\end{align*}
and $M_{n}^{\ast}(\{\mathbf{r}:m_{n}^{\ast}(\mathbf{r})\leq m_{n}^{\ast}%
(T_{n}(x)/n)\}\cap H_{n}^{c}\cap\Theta_{k,c})\leq M_{n}^{\ast}(H_{n}^{c}%
\cap\Theta_{k,c}).$ Using Lemma \ref{lem:mnstar-approx} and the bound
$|\pi(\mathbf{\theta})-\pi(\mathbf{r})|\leq B,$ when $n>(1+Bk!)^{2},$ then%
\[
M_{n}^{\ast}(H_{n}^{c}\cap\Theta_{k,c})\leq\frac{n!n^{k}}{(n+k)!}\Pi(H_{n}%
^{c}\cap\Theta_{k,c})+\frac{1}{n^{s}}+B\text{vol}(H_{n}^{c}\cap\Theta_{k,c}).
\]
So $M_{n}^{\ast}(H_{n}^{c}\cap\Theta_{k,c})\rightarrow0$ and $M_{n}^{\ast
}(\{\mathbf{r}:m_{n}^{\ast}(\mathbf{r})\leq m_{n}^{\ast}(T_{n}(x)/n)\}\cap
H_{n}^{c}\cap\Theta_{k,c})\rightarrow0$ as $n\rightarrow\infty$.

Now put
\[
A_{n}=M_{n}^{\ast}(\{\mathbf{r}\in H_{n}:m_{n}^{\ast}(\mathbf{r})\leq
m_{n}^{\ast}(T_{n}(x)/n)\}).
\]
There exist $N_{4}>N_{3}$ such that for all $n>N_{4},$ then $0<(1-n!n^{k}%
/(n+k)!)B\leq\epsilon/6.$ Then, for all $n>N_{4},$ by Lemma
\ref{lem:mnstar-approx} and the definition of $H_{n},$ for all $\mathbf{r}\in
H_{n}$,
\begin{align}
m_{n}^{\ast}(\mathbf{r}) &  \geq\frac{n!n^{k}}{(n+k)!}\pi(\mathbf{r})-\frac
{1}{n^{s}}-\sup_{\mathbf{\theta}\in G_{n,s}(\mathbf{r})\cap\Theta_{k,c}}%
|\pi(\mathbf{\theta})-\pi(\mathbf{r})|\nonumber\\
&  =\pi(\mathbf{r})-\left(  1-\frac{n!n^{k}}{(n+k)!}\right)  \pi
(\mathbf{r})-\frac{1}{n^{s}}-\sup_{\mathbf{\theta}\in G_{n,s}(\mathbf{r}%
)\cap\Theta_{k,c}}|\pi(\mathbf{\theta})-\pi(\mathbf{r})|\nonumber\\
&  \geq\pi(\mathbf{r})-\epsilon/2.\label{ineq2}%
\end{align}
By (\ref{ineq1}) and (\ref{ineq2}), for all $n>N_{4},$
\[
A_{n}\leq M_{n}^{\ast}(\{\mathbf{r}\in H_{n}:\pi(\mathbf{r})\leq\pi
(\theta^{\text{true}})+\epsilon\})\leq M_{n}^{\ast}(\{\mathbf{r}%
:\pi(\mathbf{r})\leq\pi(\theta^{\text{true}})+\epsilon\}).
\]

Putting $C_{\epsilon}(\theta^{\text{true}})=\{\mathbf{r}:\pi(\mathbf{r}%
)\leq\pi(\theta^{\text{true}})+\epsilon\},$ we have%
\begin{align*}
&  \left\vert M_{n}^{\ast}(C_{\epsilon}(\theta^{\text{true}}))-\Pi
(C_{\epsilon}(\theta^{\text{true}})\right\vert =\left\vert \int_{C_{\epsilon
}(\theta^{\text{true}})}\left(  m_{n}^{\ast}(\mathbf{\theta})-\pi
(\mathbf{\theta})\right)  \,d\mathbf{\theta}\right\vert \\
&  \leq\int_{R^{k}}\left\vert m_{n}^{\ast}(\mathbf{\theta})-\pi(\mathbf{\theta
})\right\vert \,d\mathbf{\theta\rightarrow0}%
\end{align*}
as $n\rightarrow\infty$ by Corollary \ref{cor:mnstar-approx} and Scheff\'{e}'s
theorem. Therefore, $\limsup_{n\rightarrow\infty}A_{n}\leq\Pi(\pi
(\mathbf{r})\leq\pi(\theta^{\text{true}})+\epsilon).$ Similarly, a lower bound
is obtained as $\liminf_{n\rightarrow\infty}A_{n}\geq\Pi(\pi(\mathbf{r}%
)\leq\pi(\theta^{\text{true}})-\epsilon).$ Therefore,%
\[
\limsup_{n\rightarrow\infty}A_{n}-\liminf_{n\rightarrow\infty}A_{n}\leq
\Pi(|\pi(\mathbf{\theta})-\pi(\theta^{\text{true}})|\leq\epsilon)\leq\eta.
\]
Since $\eta>0$ is arbitrary, $\lim_{n\rightarrow\infty}A_{n}$ $=\Pi
(\{\mathbf{r}\in\Theta_{k,c}:\pi(\mathbf{r})\leq\pi(\theta^{\text{true}})\})$
and the proof is complete.

\section*{References}

\noindent Al-Labadi, L. and Evans, M. (2016) Prior-based model checking.
arXiv:1606.08106 and to appear in the Canadian Journal of
Statistics.\smallskip

\noindent Al-Labadi, L., Baskurt, Z and Evans, M. (2017) Goodness of fit for
the logistic regression model using relative belief. Journal of Statistical
Distributions and Applications, 4:17.\smallskip

\noindent Briegel, H.-J., Englert, B.-G., Sterpi, N. and Walther, H. (1994)
One-atom maser: Statistics of detector clicks. Physical Review A, 49, 4,
2962-2985.\smallskip\ 

\noindent Chuang-Stein, C. and Agresti, A. (1997) A review of tests for
detecting a monotone dose--response relationship with ordinal response data.
Statistics in Medicine, 16, 22, 2599 - 2618.\smallskip

\noindent Devroye, L. (1987) A Course in Density Estimation.
Birkhauser.\smallskip

\noindent Evans, M. (2015) Measuring Statistical Evidence Using Relative
Belief. Monographs on Statistics and Applied Probability 144, CRC Press,
Taylor \& Francis Group.\smallskip

\noindent Evans, M. and Jang, G-H. (2011a). Weak informativity and the
information in one prior relative to another. Statistical Science, 26, 3,
423-439.\smallskip

\noindent Evans, M. and Jang, G-H. (2011b) A limit result for the prior
predictive applied to checking for prior-data conflict. Statistics and
Probability Letters, 81, 1034-1038.\smallskip

\noindent Evans, M. and Moshonov, H. (2006) Checking for prior-data conflict.
Bayesian Analysis, Volume 1, Number 4, 893-914.\smallskip

\noindent Evans, M. and Moshonov, H., (2007) Checking for prior-data conflict
with hierarchically specified priors. Bayesian Statistics and its
Applications, eds. A.K. Upadhyay, U. Singh, D. Dey, Anamaya Publishers, New
Delhi, 145-159.\smallskip

\noindent Evans, M. and Swartz, T. (1998) Random variable generation using
concavity properties of transformed densities. J. of Computational and
Graphical Statistics, 7, 4, 514-528.\smallskip

\noindent Izs\'{a}k, F. (2006) Maximum likelihood estimation for constrained
parameters of multinomial distributions---Application to Zipf--Mandelbrot
models. Computational Statistics and Data Analysis , 51, 3,
1575-1583.\smallskip

\noindent Len, Y.\ L., Dai, J., Englert, B-G and Krivitsky, L. A. Unambiguous
path discrimination in a two-path interferometer. arXiv:
1708.01408v2.\smallskip

\noindent Papp, L. (1992) Drosophilid assemblages in mountain creek valleys in
Hungary (Diptera: Drosophilidae) I. Folia Entomol. Hung. 53,
139--153.\smallskip

\noindent Robbins, H. (1955), A remark on Stirling's formula. The American
Mathematical Monthly, 62, 26-29.\smallskip

\noindent Shang, J., Ng, H. K., Sehrawat, A., Li, X. and Englert, B.-G. (2013)
Optimal error regions for quantum state estimation. New Journal of Physics.
15, 123026.

\end{document}